\documentclass[11pt,leqno]{amsart}
\usepackage{amsmath,amssymb,amsthm,graphicx}

\begin{document}
\title[Trigonal minimal surfaces in flat tori]
{Trigonal minimal surfaces in flat tori
\footnote{
$\quad \> \,$Key Words and Phrases: Minimal surfaces, flat tori.\\
}}
\author{Toshihiro Shoda}
\address{Faculty of Culture and Education, 
Saga University, 
1 Honjo-cho, Saga-city, Saga, 840-8502, Japan}
\email{tshoda@cc.saga-u.ac.jp}
\maketitle
\begin{abstract}
In this paper, we study trigonal minimal surfaces in flat tori. 
First, we show a topological obstruction similar to that of hyperelliptic minimal surfaces. 
Actually, the genus of trigonal minimal surface in $3$-dimensional flat torus must be 
$1$ (mod $3$). 
Next, we construct an explicit example in the higher codimensional case. 
This surface satisfies good properties and is theoretically distinct from traditional examples. 
\end{abstract}
\newtheorem{thm}{Theorem}[section]
\newtheorem{prop}{Proposition}[section]
\newtheorem{cor}{Corollary}[section]
\newtheorem{lem}{Lemma}[section]
\newtheorem{rem}{Remark}[section]
\newtheorem{main}{Main Theorem}
\newtheorem*{ex}{Example}
\newtheorem*{first_{10}}{$\dfrac{1-z^6}{w^2}\,dz$-case}
\newtheorem*{second_{10}}{$\dfrac{i\,(1+z^6)}{w^2}\,dz$-case}
\newtheorem*{third_{10}}{$\dfrac{z^5+z}{w^2}$-case}
\newtheorem*{firth_{10}}{$\dfrac{i\,(z^5-z)}{w^2}$-case}

\section{Introduction}

Let $f:M_g\longrightarrow \textbf{R}^n/\Lambda$ be a minimal immersion of a compact surface 
with genus $g$ into an $n$-dimensional flat torus which does not lie in any subtorus of 
$\textbf{R}^n/\Lambda$ (it can be replaced by an $n$-periodic minimal immersion from 
some covering space of $M_g$ into $\textbf{R}^n$). 
The conformal structure induced by the isothermal coordinates makes $M_g$ a Riemann surface 
and $f$ is called a conformal minimal immersion. 
Before we refer to the back grounds, we state the following result (p.884 \cite{Me}):
\begin{thm}[Generalized  Weierstrass Representation]$\,$

If $f:M_g\longrightarrow\textbf{R}^n/\Lambda$ is a conformal minimal immersion then, 
after a translation, $f$ can be represented by
\begin{equation*}
f(p)=\Re\int_{p_0}^p (\omega_1, \omega_2, \dots, \omega_n)^T \quad Mod \>\> \Lambda,
\end{equation*}
where $p_0$ is a fixed point in $M_g$, superscript $T$ means the transposed matrix, and 
$\omega_1$, $\omega_2$, $\dots$, $\omega_n$ are holomorphic differentials on $M_g$ satisfying
\begin{align}
&\omega_1, \omega_2, \dots, \omega_n{\rm\,\, have\,\, no\,\, common\,\, zeros},\\
&\sum_{k=1}^n \omega_k^2=0 \,,\\
&\left\{ \Re \int_{\gamma} (\omega_1, \omega_2, \dots, \omega_n)^T \,\,|\,\, \gamma \in 
H_1 (M_g, \textbf{Z}) \right\}{\rm is\,\, a\,\, sublattice\,\, of\,\, }\Lambda.
\end{align}
Conversely, if $\omega_1$, $\omega_2$, $\cdots$, $\omega_n$ 
are holomorphic differentials satisfying $(1)$, $(2)$, and $(3)$, 
then $f$, defined as above, is a conformal minimal immersion. 
\end{thm}
Condition (3) is called the period condition and guarantees the well-definedness of 
the path integral. 

Using the notations as above, we define the associate surface $f_{\theta}$: 
\[f_{\theta}(p):=\Re\int_{p_0}^p e^{i\,\theta}\,(\omega_1, \omega_2, \dots, \omega_n)^T.\]
In particular, the associate surface $f_{\pi/2}$ is called the conjugate surface. 
We note that $f_{\theta}$ may not be well-defined for any torus, 
even though $f=f_0$ is well-defined on $M_g$. 
Recall that the Gauss map $G$ is given as a holomorphic map from $M_g$ to the quadric 
$Q_{n-2}:=\{ [w_1,\,w_2,\,\cdots,\,w_n] \in \textbf{C}P^{n-1} | \sum_{k=1}^n (w_k)^2 = 0 \}$:
\begin{align*}
G:\>&M_g\longrightarrow Q_{n-2}\subset \textbf{C}P^{n-1}\\
&\>\>p\longmapsto [\omega_1(p),\,\omega_2(p),\,\cdots,\,\omega_n(p)]
\end{align*}

One of the beautiful classical theorems on compact Riemann surfaces states that every compact 
Riemann surface of positive genus is holomorphically embedded in the Jacobian by 
the Abel-Jacobi map: take a basis $\{\eta_1,\eta_2,\cdots,\eta_g\}$ of the space of holomorphic 
differentials of $M_g$, and consider 
\[\Lambda_{\eta}=\left\{\Re \int_{\gamma}\left(\eta_1,\cdots,\eta_g,-i\,\eta_1,\cdots,-i\,\eta_g 
\right)^T |\>\gamma\in H_1(M_g,\textbf{Z})\right\}.\]
The Jacobian ${\rm Jac}(M_g)$ is the complex torus 
represented by $\textbf{C}^g/\Lambda_{\eta}$ and the holomorphic 
embedding $j:M_g\longrightarrow {\rm Jac}(M_g)$ defined by 
\[j(p)=\Re\int^p_{p_0}\left(\eta_1,\cdots,\eta_g,-i\,\eta_1,\cdots,-i\,\eta_g\right)^T\]
is called the Abel-Jacobi map. 
The Jacobian satisfies the following universal property (p.5 \cite{N-S2}): 

\begin{thm}$\,$

Given $f:M_g\longrightarrow \textbf{R}^n/\Lambda$ a conformal 
minimal immersion, we may assume $f(p_0)=0$. Then there exists 
a real homomorphism $h$ from ${\rm Jac}(M_g)$ 
to $\textbf{R}^n/\Lambda$ so that $f=h\circ j${\rm :}

\begin{picture}(300,45)
 \put(80,-7){\vector(1,0){170}}
 \put(50,-10){$M_g$}
 \put(270,-10){$\textbf{R}^n/\Lambda$}
 \put(150,40){{\rm Jac}$(M_g)$}
 \put(70,0){\vector(2,1){70}}
 \put(200,35){\vector(2,-1){70}}
 \put(100,25){$j$}
 \put(230,25){$\exists\, h$}
 \put(165,-18){$f$}
 \put(165,15){$\circlearrowright$}
\end{picture}
\\
\end{thm}

The Abel-Jacobi map plays an important role in the theory of algebraic curves 
(e.g. Torelli's Theorem, Schottoky problem, etc.) and so, by Theorem 1.2, 
it is useful to study minimal surfaces from the point of 
view of the theory of algebraic curves. In fact, Ejiri \cite{E} has translated Schottoky problem 
into a differential geometric situation and has studied the Moduli space of compact minimal 
surfaces in flat tori. 

Algebraic curves can be divided into non-hyperelliptic curves and hyperelliptic curves, and 
in particular, there is a topological obstruction to hyperelliptic minimal surfaces in 
$3$-dimensional flat tori. 
Actually, a hyperelliptic minimal surface of even genus cannot be minimally immersed into any 
$3$-dimensional flat torus (Theorem 3.3 in \cite{Me}). 

Note that every compact Riemann surface can be represented as a branched $d$-cover of the sphere 
for some $d\geq 1$,  and it is reasonable to ask whether there is a topological obstruction or not 
for $d> 2$ (a Riemann surface with $d=2$ is a hyperelliptic curve, and hence we omit 
the case $d=2$). 
Now we consider this problem for $d=3$. 
Recall that a non-hyperelliptic curve with $d=3$ is called a trigonal curve. 

Our first result is the following topological obstruction:  
\begin{main}$\,$

Let $f:M_g \longrightarrow \textbf{R}^3/\Lambda$ be a conformal minimal immersion of a 
trigonal Riemann surface $M_g$ with genus $g$. Then, $g=3r+1$ holds for some $r\geq 1$. 
Therefore, a trigonal Riemann surface of genus $0$ or $2$ {\rm (}mod $3${\rm )} cannot be 
minimally immersed into any $3$-dimensional flat torus.  
\end{main}
\begin{rem}$\,$

There are some examples of trigonal minimal surfaces in $3$-dimensional flat tori{\rm :} 
{\rm (i)} A previous example \cite{Shoda} {\rm (}$r=1${\rm )}, 
{\rm (ii)} Schoen's I-WP surface \cite{K}, \cite{A} {\rm (}$r=3${\rm )}, where $r$ 
is given in Main Theorem 1. 
\end{rem}

Next, we consider the higher codimensional case. In 1976, Nagano-Smyth \cite{N-S3} 
constructed compact minimal surfaces in $n$-dimensional flat tori abstractly. 
On the other hand, there are few explicit examples, and thus we will construct an example of 
trigonal minimal surfaces in the simplest case, that is, in case $n=4$. 
Actually, an example of trigonal minimal surfaces of genus $10$ in $4$-dimensional flat tori 
is given (see 3.2). 
In general, the most difficult part in constructing examples comes from the period condition. 
It is not always possible to solve the period condition, and we cannot even calculate 
the periods. We overcome this problem through the following process: 
(i) taking a suitable Riemann surface 
and considering its symmetries (see 3.1), 
(ii) calculating the periods by the functional theoretical 
techniques (see Lemma 3.1), and (iii) finding a relation between one period and another one 
(see Lemma 3.2). 

Now, we consider our example in view of Nagano-Smyth's discussion \cite{N-S2}. 
In fact, this example satisfies the following properties: 
(a) the conjugate surface $f_{\pi/2}:M_{10}\longrightarrow \textbf{R}^4/\Lambda_{\pi/2}$ 
is well-defined, 
(b) the associate surfaces $f_{\theta}:M_{10}\longrightarrow \textbf{R}^4/\Lambda_{\theta}$ 
are also well-defined for a countable dense set of angles $e^{i\,\theta}\subset S^1$, 
(c) this example is homologous to $0$ in the $4$-torus. 
\begin{rem}[section 6 in \cite{A-Pi}]$\,$

Arezzo-Pirola showed the existence of minimal surfaces which are not homologous to $0$ in the tori. 
\end{rem}
Here, we review the Nagano-Smyth's argument. Let $S_f(M_g)$ be a subgroup of the automorphism 
group of $M_g$, and say $f$ has symmetry $S_f(M_g)$ if and only if 
$S_f(M_g)$ extends under $f$ 
to a group of affine transformations of $\textbf{R}^n/\Lambda$. When the corresponding 
linear representation of $S_f(M_g)$ is irreducible, we say $f$ has 
irreducible symmetry $S_f(M_g)$. 
If the complexification of this representation is also irreducible, then we say $f$ has 
absolutely irreducible symmetry $S_f(M_g)$. 
Nagano-Smyth showed that (i) if $f$ has absolutely irreducible symmetry, 
then $f$ satisfies (c) (Theorem 2 in \cite{N-S2}), (ii) if we assume the above irreducible 
conditions by Weyl group and some conditions, then $f$ satisfies (a) and (b) 
(Theorem 5 in \cite{N-S2}). 
However, our example satisfies (a), (b), (c), and has only reducible symmetry (see section 3). 
Therefore, the converse propositions of Nagano-Smyth's results do not hold respectively. 

Finally, we state the above results as follows: 
\begin{main}$\,$

There exists a trigonal minimal surface of genus $10$ 
in a $4$-dimensional flat torus satisfying the following properties{\rm:} 
{\rm (i)} the conjugate surface $f_{\pi/2}:M_{10}\longrightarrow \textbf{R}^4/\Lambda_{\pi/2}$ 
is well-defined, {\rm (ii)} the associate surfaces $f_{\theta}:M_{10}\longrightarrow 
\textbf{R}^4/\Lambda_{\theta}$ are well-defined for a countable dense set of angles 
$e^{i\,\theta}\subset S^1$, 
{\rm (iii)} this surface is homologous to $0$ in the torus and has only reducible symmetry. 
\end{main}

\section{A topological obstruction to trigonal minimal surfaces in $3$-dimensional flat tori}

In this section, we prove Main Theorem 1. First, we review the spinor representation of minimal 
surfaces \cite{K-S}. 
Let $f:M_g\longrightarrow \textbf{R}^3/\Lambda$ be a conformal minimal immersion defined 
by $f(p)=\int^p_{p_0}(\omega_1,\omega_2,\omega_3)^T$. Then, the Gauss map $G$ commutes 
the following diagram: 

\begin{picture}(300,50)
 \put(90,33){\vector(1,0){150}}
 \put(70,30){$M_g$}
 \put(245,30){$Q_1\subset \textbf{C}P^2$}
 \put(155,-10){$\textbf{C}P^1$}
 \put(90,25){\vector(2,-1){60}}
 \put(195,8){$\cong$}
 \put(180,-5){\vector(2,1){60}}
 \put(105,5){$\varphi$}
 \put(220,5){$V$}
 \put(160,40){$G$}
 \put(160,12){$\circlearrowright$}
\end{picture}
$\,$\\
\\
where $V$ is the Veronese embedding given by 
\[V(s_1,s_2)=(s_1^2-s_2^2,i\,(s_1^2+s_2^2),2\,s_1s_2)^T.\]
Let $\mathcal{O}_{\textbf{C}P^n}(1)$ be the hyperplane bundle on $\textbf{C}P^n$ and 
$L=\varphi^*(\mathcal{O}_{\textbf{C}P^1}(1))$ its pull-back to $M_g$. Note that 
$L^2=G^*(\mathcal{O}_{\textbf{C}P^2}(1))$ is the canonical bundle $K$ of $M_g$, i.e. $L$ defines 
a spin structure on $M_g$. Besides, again by pull-back, we find two holomorphic sections 
$t_1$, $t_2$ of $L$ which have no common zeros such that 
\[(\omega_1,\omega_2,\omega_3)^T=(t_1^2-t_2^2,i\,(t_1^2+t_2^2),2\,t_1t_2)^T.\]
The meromorphic function $t_2/t_1$ can be identified with the usual Gauss map 
$M_g\longrightarrow S^2\cong \textbf{C}\cup \{\infty\}$ (see \cite{H-Os}). 

Next, we give some notations for a linear series (or system). The standard terminology is as 
in \cite{A-C-G-H}. 
Given a divisor $D$ on $M_g$, the complete linear series $|D|$ is the set of effective 
divisors linearly equivalent to $D$. We then have an identification between $|D|$ and 
projectivization of the space of holomorphic sections of the line bundle defined by $D$: 
$|D|\cong PH^0(M_g,\mathcal{O}(D))$, and thus a complete linear series is a projective space. 
More generally, every linear subspace of a complete linear series is called a linear series. 
A linear series $PW$, where $W$ is a vector subspace of $|D|$, is said to be a 
$g^r_d$ if and only if $\deg D=d$ and $\dim W=r+1$. A $g^1_d$ is called a pencil, a $g^2_d$ a net, 
and a $g^3_d$ a web. We write $t\,g^r_d$ (resp. $|K-g^r_d|$)
for the complete linear series $|\,t\,E|$ (resp. $|K-E|$), where 
$E\in g^r_d$. By a base point of a linear series $PW$, we mean a point common to all 
divisors  of $PW$. If there are none, then we say that the linear series is base-point-free. 

Using the aboves, we prove Main Theorem 1 as follows: 
\begin{proof}$\,$

We first observe that the degree of the Gauss map is $g-1$ (Theorem 3.1 in \cite{Me}). 
It follows that $M_g$ is not trigonal in case $g=0,1,2,3$, and hence $g\geq 4$ holds. 

We may omit the case $g=4$ because it corresponds to the case $r=1$. 
Note that $M_g$ is trigonal if and only if there is a base-point-free pencil $g^1_3$ on $M_g$. 
In case $g>4$, the $g^1_3$ on $M_g$ is unique complete linear series (p.124 \cite{Sho}). 

Let $L$ be the spin bundle given by the spinor representation of $f$ and $D_L$ the 
divisor defined by $L$. Note that (1) in Theorem 1.1 implies that $|D_L|$ is base-point-free 
complete linear series. It is known that every base-point-free complete linear series $g^r_d$ 
is represented by 
\[g^r_d=r\,g^1_3\>\>\>{\rm or}\>\>\> |K-g^r_d|=r'\,g^1_3+F\>\>\> (r':=g-d+r-1), \]
where F is an effective divisor and consists of the base points of $|K-g^r_d|$ (Remark 1.2 in 
\cite{C-M}, (1.2.7) in \cite{C-K-M}). 
Applying the above fact to $|D_L|$ ($=|K-D_L|$), we obtain $|D_L|= r\,g^1_3$ for some $r>1$. 
In particular, $\deg D_L=g-1=3r=\deg (r\,g^1_3)$ follows and the proof is complete.  
\end{proof}

\section{An example of trigonal minimal surfaces in $4$-dimensional flat tori}

In this section, we prove Main Theorem 2 through the following four subsections. 
In 3.1, we explain how to determine our example. In 3.2, we calculate periods of the surface 
and solve the period condition directly. Moreover, we consider the conjugate surface and 
the associate surfaces. 
In 3.3, we prove that the surface is homologous to $0$ in the torus. 
Finally, we show that the surface has only reducible symmetry in 3.4. \\

\subsection{Construction}$\,$

In \cite{Shoda}, we constructed a trigonal minimal surface defined by $w^3=z^6-1$. We now 
consider the higher genus version of it. Let $M_g$ be the cyclic covering of a line 
(p.73 \cite{Mir}) given by 
\begin{equation*}
w^3=z^{g+2}-1\>\>\>(g=3r+1,\>\>r=1,2,3,\cdots).
\end{equation*}
Then, we can compute explicitly a basis for the space of holomorphic differentials $H^0 (M_g, K)$:
\[H^0 (M_g, K)={\rm span}\left\{
\dfrac{dz}{w^2},z\,\dfrac{dz}{w^2},\cdots,z^{2r}\,\dfrac{dz}{w^2},
\dfrac{dz}{w},z\,\dfrac{dz}{w},\cdots,z^{r-1}\,\dfrac{dz}{w}\right\}.\]
Note that $M_g$ is not well-defined in case $g=3r+2$ because Riemann-Hurwitz's formula 
does not hold, and 
$M_g$ has cusp singularity at $z=\infty$ in case $g=3r$. We make a choice the simplest case, 
that is, the case $g=3r+1$ as above. 

Now we consider the following conformal minimal immersion $f$ 
from $M_g$ into a $\textbf{R}^4/\Lambda$ ($\Lambda$ is defined later):
\[f(p)=\Re \int_{p_0}^{p} \underbrace{\left(\dfrac{1-z^{2r}}{w^2},\,
\dfrac{i\,(1+z^{2r})}{w^2},\,\dfrac{z^{2r-1}+z}{w^2},\,
i\,\dfrac{z^{2r-1}-z}{w^2}\right)^T\,dz}_{\Psi}.\]

\begin{rem}$\,$

The author found a component of the Moduli space which is corresponding to trigonal minimal 
surfaces in $4$-tori \cite{Shoda}. 
The minimal immersion $f$ defined as above is an element of the component of the Moduli 
space. 
\end{rem} 
To find symmetries of $f$, we now consider the automorphism $\varphi$ defined by 
$\varphi(z,w)=(e^{\frac{2\pi}{3(r+1)}\,i} z, w)$. 
Then, we get 
\begin{align*}
&\varphi^* \Psi=e^{\frac{2\,\pi}{3}\,i}
\begin{pmatrix}
R\left(\dfrac{2\pi}{3}-\dfrac{2\pi}{3(r+1)}\right) & 0\\
0 & R\left(\dfrac{4\pi}{3(r+1)}-\dfrac{2\pi}{3}\right) 
\end{pmatrix}\Psi,
\end{align*}
where
\begin{align*}
R(\theta):=
\begin{pmatrix}
\cos\theta&-\sin\theta\\
\sin\theta&\cos\theta
\end{pmatrix}.
\end{align*}

\begin{rem}$\,$

In case $r=1,3,7,9$, 
$\sin\left(\frac{2\pi}{3}-\frac{2\pi}{3(r+1)}\right)$ can be calculated as the following 
explicit values{\rm :} 
$\left\{\frac{\sqrt{3}}{2}, 1, \frac{\sqrt{3}+1}{2\sqrt{2}}, 
\frac{1}{2}\,\sqrt{\frac{5+\sqrt{5}}{2}}\right\}$. 
We use $\varphi^*$ to calculate periods, and we cannot solve the period condition if 
$\sin\left(\frac{2\pi}{3}-\frac{2\pi}{3(r+1)}\right)$ is too complicated. 
So, the period condition may not be solved in case $r \neq 1,3$.
These suggest that $r=1,\,3$ may be the best case to construct examples explicitly. 
\end{rem}

We now consider the case $g=10$, that is, $r=3$. 
Then $M_{10}$ is defined by $w^3=z^{12}-1$ and $f$ can be reduced to 
\[f(p)=\Re \int_{p_0}^{p} \underbrace{\left(
\dfrac{1-z^6}{w^2},\,
\dfrac{i\,(1+z^6)}{w^2},\,
\dfrac{z^5+z}{w^2},\,
\dfrac{i\,(z^5-z)}{w^2}
\right)^T\,dz}_{\Psi}. \]
$\varphi$ is given by 
$\varphi(z,w)=(e^{\frac{\pi}{6}\,i} z, w)$ 
and 
\begin{align*}
&\varphi^* \Psi
=e^{\frac{2\,\pi}{3}\,i}
\begin{pmatrix}
R\left(\dfrac{\pi}{2}\right) & 0\\
0 & R\left(-\dfrac{\pi}{3}\right) 
\end{pmatrix}
\Psi
=e^{\frac{2\,\pi}{3}\,i}
\begin{pmatrix}
0 & -1& 0 & 0 \\
1 & 0 & 0 & 0 \\
0 & 0 & \dfrac{1}{2} & \dfrac{\sqrt{3}}{2}\\
0 & 0 & -\dfrac{\sqrt{3}}{2}& \dfrac{1}{2}
\end{pmatrix}\Psi.
\end{align*}

In the next subsection, we prove the well-definedness of the following conformal minimal immersion 
of $M_{10}$ into $\textbf{R}^4/\Lambda$: 
\begin{align}
f:&\>M_{10}\longrightarrow\textbf{R}^4/\Lambda\\
\nonumber  &\>\>p\>\>\longmapsto\Re \int_{p_0}^{p} \left(
\dfrac{1-z^6}{w^2},\,
\dfrac{i\,(1+z^6)}{w^2},\,
\dfrac{z^5+z}{w^2},\,
\dfrac{i\,(z^5-z)}{w^2}
\right)^Tdz\\
\nonumber&\qquad\qquad\qquad\qquad\qquad\qquad\qquad\qquad\qquad\qquad(w^3=z^{12}-1),
\end{align}
where $\Lambda$ is given by the Beta function B(a,b): 
\[\Lambda=\begin{pmatrix}
3\,\alpha & \dfrac{3}{2}\,\alpha & 0 & 0\\
0 & \dfrac{3}{2}\,\alpha & 0 & 0\\
0 & 0 & 3\,\gamma & \dfrac{3}{2}\,\gamma\\
0 & 0 & 0 & \dfrac{\sqrt{3}}{2}\,\gamma
\end{pmatrix},\>\>
\begin{cases}
\alpha=\dfrac{1}{6\sqrt[3]{2}}\,B(2/3,1/6), \\
\gamma= \dfrac{1}{12}\,B(1/3, 1/6)\end{cases}\]

\subsection{Periods}$\,$

In this subsection, we calculate the periods of $f$ defined by (4), and consider the conjugate 
surface and the associate surfaces. 
$1$-cycles $A_1$, $A_2$, $\cdots$, $A_{10}$, $B_1$, $B_2$, $\cdots$, $B_{10}$ 
are obtained as follows: 
\begin{align*}
A_1&=\left\{(z,\,w)=(e^{i\,t},\,w(t))\,|\,t\in \left[0,\dfrac{\pi}{6}\right],\,
w\left(\dfrac{\pi}{12}\right)=-\sqrt[3]{2} \right\}\\
&\cup \left\{(z,\,w)=(e^{-i\,t},\,w(t))\,|\,
t\in \left[-\dfrac{\pi}{6},\,0\right],\,w\left(-\dfrac{\pi}{12}\right)=-\sqrt[3]{2}\,
e^{\frac{2\,\pi}{3}\,i}\right\},\\
A_2&=\left\{(z,\,w)=(e^{i\,t},\,w(t))\,|\,t\in \left[0,\,\dfrac{\pi}{6}\right],
\,w\left(\dfrac{\pi}{12}\right)=-\sqrt[3]{2} \right\}\\
&\cup \left\{(z,\,w)=(e^{-i\,t},\,w(t))\,|\,
t\in \left[-\dfrac{\pi}{6},\,0\right],\,w\left(-\dfrac{\pi}{12}\right)=-\sqrt[3]{2}\,
(e^{\frac{2\,\pi}{3}\,i})^2\right\},
\end{align*}
and 
\begin{align*}
A_{2k+1}=\varphi^2(A_{2k-1})&, \>\>A_{2k+2}=\varphi^2(A_{2k})\qquad (1\leq k\leq 4),\\
B_{2l-1}=\varphi(A_{2l-1})&, \>\>B_{2l}=\varphi(A_{2l})\qquad (1\leq l\leq 5).  
\end{align*}

First, we calculate periods along $A_1$ and $A_2$: 
\begin{lem}
\begin{align*}
\begin{pmatrix}
\int_{A_1} \dfrac{1-z^6}{w^2}\,dz\\
\int_{A_1} \dfrac{i\,(1+z^6)}{w^2}\,dz\\
\int_{A_1} \dfrac{z^5+z}{w^2}\,dz\\
\int_{A_1} \dfrac{i\,(z^5-z)}{w^2}\,dz
\end{pmatrix}
&=\begin{pmatrix}
\dfrac{1+e^{\frac{\pi}{3}\,i}}{6\sqrt[3]{2}}\,B(2/3,1/6)\\
-\dfrac{1+e^{\frac{\pi}{3}\,i}}{6\sqrt[3]{2}}B(2/3,1/6) \\
\dfrac{i\,(1+e^{\frac{\pi}{3}\,i})}{4\sqrt{3}}\, B(1/3,1/6)\\
-i\,(1+e^{\frac{\pi}{3}\,i})\,\int^1_{\frac{1}{2}}\dfrac{dt}
{\sqrt[3]{4\,(1-t^2)(4t^2-1)^2}}
\end{pmatrix},\\
\begin{pmatrix}
\int_{A_2} \dfrac{1-z^6}{w^2}\,dz\\
\int_{A_2} \dfrac{i\,(1+z^6)}{w^2}\,dz\\
\int_{A_2} \dfrac{z^5+z}{w^2}\,dz\\
\int_{A_2} \dfrac{i\,(z^5-z)}{w^2}\,dz
\end{pmatrix}
&=\begin{pmatrix}
\dfrac{e^{\frac{2\pi}{3}\,i}+e^{\frac{\pi}{3}\,i}}{6\sqrt[3]{2}}\,B(2/3,1/6)\\ 
-\dfrac{e^{\frac{2\pi}{3}\,i}+e^{\frac{\pi}{3}\,i}}{6\sqrt[3]{2}}B(2/3,1/6)\\
\dfrac{i\,(e^{\frac{\pi}{3}\,i}+e^{\frac{2\pi}{3}\,i})}{4\sqrt{3}}\,B(1/3,1/6)\\
-i\,(e^{\frac{\pi}{3}\,i}+e^{\frac{2\pi}{3}\,i})\,\int^1_{\frac{1}{2}}\dfrac{dt}
{\sqrt[3]{4\,(1-t^2)(4t^2-1)^2}}
\end{pmatrix}.
\end{align*}
\end{lem}

\begin{proof}

\begin{first_{10}}
\end{first_{10}}
First, we consider 
$A_1$-period. In case $t: 0\longmapsto \pi/6$, we set 
\[\eta=\dfrac{1}{2}\left(z^3+\dfrac{1}{z^3}\right)=
\dfrac{1}{2}\left(e^{3i\,t}+e^{-3i\,t}\right)=\cos (3t).\] 
Then, $d\eta=3\,\dfrac{z^6-1}{2\,z^4}\,dz$, and thereby 
\[\dfrac{1-z^6}{w^2}\,dz=\dfrac{1-z^6}{w^2}\,\dfrac{2\,z^4}{3\,(z^6-1)}\,d\eta
=-\dfrac{2}{3}\,\dfrac{z^4}{w^2}\,d\eta.\]  
To calculate $z^4/w^2$, we consider 
\[\left(\dfrac{z^4}{w^2}\right)^3=\dfrac{z^{12}}{(z^{12}-1)^2}=
\dfrac{1}{z^{12}+\frac{1}{z^{12}}-2}. \]
Note that 
$z^6+\dfrac{1}{z^6}=4\,\eta^2-2$ and thus 
\[z^{12}+\dfrac{1}{z^{12}}=(4\,\eta^2-2)^2-2
=16\eta^4-16\eta^2+2.\] So, we obtain 
\[\left(\dfrac{z^4}{w^2}\right)^3=-\dfrac{1}{16\,\eta^2\,(1-\eta^2)}<0. \]
Here, we take a suitable branch for $z^4/w^2$. 
\[\dfrac{z^4}{w^2}\left(\dfrac{\pi}{12}\right)=\dfrac{e^{\frac{\pi}{3}\,i}}
{(-\sqrt[3]{2})^2}=\dfrac{e^{\frac{\pi}{3}\,i}}{\sqrt[3]{4}}. \]
Choosing a branch $\sqrt[3]{\eta^2(1-\eta^2)}>0$, we get 
\[\dfrac{z^4}{w^2}=\dfrac{e^{\frac{\pi}{3}\,i}}{2\sqrt[3]{2}\,
\sqrt[3]{\eta^2(1-\eta^2)}}.\] It follows that 
\[\dfrac{1-z^6}{w^2}\,dz=-\dfrac{1}{3\sqrt[3]{2}}\,
\dfrac{e^{\frac{\pi}{3}\,i}}{ \sqrt[3]{\eta^2(1-\eta^2)}}\,d\eta. \]
In case $t:-\pi/6\longmapsto 0$, by the similar arguments, we obtain 
\[\dfrac{1-z^6}{w^2}\,dz=\dfrac{1}{3\sqrt[3]{2}}\,
\dfrac{1}{\sqrt[3]{\eta^2(1-\eta^2)}}\,d\eta. \]
Therefore, we obtain 
\begin{align*}
\int_{A_1} \dfrac{1-z^6}{w^2}\,dz&=\int^0_1-\dfrac{1}{3\sqrt[3]{2}}\,
\dfrac{e^{\frac{\pi}{3}\,i}}{\sqrt[3]{\eta^2(1-\eta^2)}}\,d\eta+
\int^1_{0}\dfrac{1}{3\sqrt[3]{2}}\,
\dfrac{1}{\sqrt[3]{\eta^2(1-\eta^2)}}\,d\eta\\
&=\dfrac{1+e^{\frac{\pi}{3}\,i}}{3\sqrt[3]{2}}\int^1_0
\dfrac{d\eta}{\sqrt[3]{\eta^2(1-\eta^2)}}
\underbrace{=}_{t=\eta^2}\dfrac{1+e^{\frac{\pi}{3}\,i}}{6\sqrt[3]{2}}\int^1_0
t^{-5/6}(1-t)^{-1/3}\,dt\\
&=\dfrac{1+e^{\frac{\pi}{3}\,i}}{6\sqrt[3]{2}}\,B(2/3,1/6).
\end{align*}
Similarly, $A_2$-period is calculated as follows: 
\begin{align*}
\int_{A_2} \dfrac{1-z^6}{w^2}\,dz&=\int^0_1-\dfrac{1}{3\sqrt[3]{2}}\,
\dfrac{e^{\frac{\pi}{3}\,i}}{\sqrt[3]{\eta^2(1-\eta^2)}}\,d\eta+
\int^1_{0}\dfrac{1}{3\sqrt[3]{2}}\,
\dfrac{e^{\frac{2\pi}{3}\,i}}{\sqrt[3]{\eta^2(1-\eta^2)}}\,d\eta\\
&=\dfrac{e^{\frac{2\pi}{3}\,i}+e^{\frac{\pi}{3}\,i}}{3\sqrt[3]{2}}\int^1_0
\dfrac{d\eta}{\sqrt[3]{\eta^2(1-\eta^2)}}
=\dfrac{e^{\frac{2\pi}{3}\,i}+e^{\frac{\pi}{3}\,i}}{6\sqrt[3]{2}}\,B(2/3,1/6). 
\end{align*}
$\,$\\
\begin{second_{10}}
\end{second_{10}}
First, we consider $A_1$-period. 
In case $t:0\longmapsto \pi/6$, we set 
\[\eta=-\dfrac{i}{2}\left(z^3-\dfrac{1}{z^3}\right)=
-\dfrac{i}{2}\left(e^{3i\,t}-e^{-3i\,t}\right)=\sin (3t).\] 
Then, $d\eta=-3i\,\dfrac{z^6+1}{2\,z^4}\,dz$, and thereby 
\[\dfrac{i\,(1+z^6)}{w^2}\,dz=\dfrac{i\,(1+z^6)}{w^2}\,\dfrac{2\,z^4}{-3i\,(z^6+1)}\,d\eta
=-\dfrac{2}{3}\,\dfrac{z^4}{w^2}\,d\eta.\] 
To calculate $z^4/w^2$, we consider
\[\left(\dfrac{z^4}{w^2}\right)^3=\dfrac{z^{12}}{(z^{12}-1)^2}=
\dfrac{1}{z^{12}+\frac{1}{z^{12}}-2}.\]
Note that 
$z^6+\dfrac{1}{z^6}=-4\,\eta^2+2$ and thus \[z^{12}+\dfrac{1}{z^{12}}=(-4\,\eta^2+2)^2-2
=16\eta^4-16\eta^2+2.\] So, we obtain 
\[\left(\dfrac{z^4}{w^2}\right)^3=-\dfrac{1}{16\,\eta^2\,(1-\eta^2)}<0.\]
Here, we take a suitable branch for $z^4/w^2$. 
\[\dfrac{z^4}{w^2}\left(\dfrac{\pi}{12}\right)=\dfrac{e^{\frac{\pi}{3}\,i}}
{(-\sqrt[3]{2})^2}=\dfrac{e^{\frac{\pi}{3}\,i}}{\sqrt[3]{4}}.\]
Choosing a branch $\sqrt[3]{\eta^2(1-\eta^2)}>0$, we get 
\[\dfrac{z^4}{w^2}=\dfrac{e^{\frac{\pi}{3}\,i}}{2\sqrt[3]{2}\,
\sqrt[3]{\eta^2(1-\eta^2)}}.\]
It follows that 
\[\dfrac{i\,(1+z^6)}{w^2}\,dz=-\dfrac{1}{3\sqrt[3]{2}}\,
\dfrac{e^{\frac{\pi}{3}\,i}}{\sqrt[3]{\eta^2(1-\eta^2)}}\,d\eta. \]
In case $t:-\pi/6\longmapsto 0$, by the similar arguments, we obtain 
\[\dfrac{i\,(1+z^6)}{w^2}\,dz=\dfrac{1}{3\sqrt[3]{2}}\,
\dfrac{1}{\sqrt[3]{\eta^2(1-\eta^2)}}\,d\eta. \]
Therefore, we obtain 
\begin{align*}
\int_{A_1} \dfrac{i\,(1+z^6)}{w^2}\,dz&=\int^1_0-\dfrac{1}{3\sqrt[3]{2}}\,
\dfrac{e^{\frac{\pi}{3}\,i}}{\sqrt[3]{\eta^2(1-\eta^2)}}\,d\eta+
\int^0_{1}\dfrac{1}{3\sqrt[3]{2}}\,
\dfrac{1}{\sqrt[3]{\eta^2(1-\eta^2)}}\,d\eta\\
&=-\dfrac{1+e^{\frac{\pi}{3}\,i}}{3\sqrt[3]{2}}\int^1_0
\dfrac{d\eta}{\sqrt[3]{\eta^2(1-\eta^2)}}
=-\dfrac{1+e^{\frac{\pi}{3}\,i}}{6\sqrt[3]{2}}B(2/3,1/6). 
\end{align*}
Similarly, $A_2$-period is calculated as follows: 
\begin{align*}
\int_{A_2} \dfrac{i\,(1+z^6)}{w^2}\,dz&=\int^1_0-\dfrac{1}{3\sqrt[3]{2}}\,
\dfrac{e^{\frac{\pi}{3}\,i}}{\sqrt[3]{\eta^2(1-\eta^2)}}\,d\eta+
\int^0_1\dfrac{1}{3\sqrt[3]{2}}\,
\dfrac{e^{\frac{2\pi}{3}\,i}}{\sqrt[3]{\eta^2(1-\eta^2)}}\,d\eta\\
&=-\dfrac{e^{\frac{2\pi}{3}\,i}+e^{\frac{\pi}{3}\,i}}{3\sqrt[3]{2}}\int^1_0
\dfrac{d\eta}{\sqrt[3]{\eta^2(1-\eta^2)}}
=-\dfrac{e^{\frac{2\pi}{3}\,i}+e^{\frac{\pi}{3}\,i}}{6\sqrt[3]{2}}B(2/3,1/6). 
\end{align*}

\begin{third_{10}}
\end{third_{10}}
First, we consider $A_1$-period. 
In case $0\longmapsto \pi/6$, we set 
\[\eta=\dfrac{1}{2\,i}\,\left(z^2-\dfrac{1}{z^2}\right)=
\sin (2t).\]
Then, $d\eta=\dfrac{1}{i}\,\dfrac{1+z^4}{z^3}\,dz$, and thereby 
\[\dfrac{z^5+z}{w^2}\,dz=\dfrac{z(1+z^4)}{w^2}\,\dfrac{i\,z^3}{1+z^4}\,d\eta
=i\,\dfrac{z^4}{w^2}\,d\eta.\] 
To calculate $z^4/w^2$, we consider 
\[\left(\dfrac{z^4}{w^2}\right)^3=\dfrac{z^{12}}{(z^{12}-1)^2}
=\dfrac{1}{z^{12}+\dfrac{1}{z^{12}}-2}. \]
Note that 
$z^2-\dfrac{1}{z^2}=2i\,\eta$, 
$z^4+\dfrac{1}{z^4}=2-4\eta^2$, and thus 
\[z^8+\dfrac{1}{z^8}=(2-4\eta^2)^2-2=2-16\eta^2+16\eta^4.\] It follows that 
\[z^{12}+\dfrac{1}{z^{12}}-2=-4\eta^2(3-4\eta^2)^2.\]
So, we obtain 
\[\left(\dfrac{z^4}{w^2}\right)^3=-\dfrac{1}{4\eta^2(3-4\eta^2)^2}<0. \]
Here, we take a suitable branch for $z^4/w^2$. 
\[\dfrac{z^4}{w^2}\left(\dfrac{\pi}{12}\right)=\dfrac{e^{\frac{\pi}{3}\,i}}
{\sqrt[3]{4}}. \]
Choosing a branch $\sqrt[3]{4\eta^2(3-4\eta^2)^2}>0$, we get 
\[\dfrac{z^4}{w^2}=\dfrac{e^{\frac{\pi}{3}\,i}}{\sqrt[3]{4\,\eta^2(3-4\eta^2)^2}}.\] 
It follows that 
\[\dfrac{z^5+z}{w^2}\,dz=i\,
\dfrac{e^{\frac{\pi}{3}\,i}}{\sqrt[3]{4\eta^2(3-4\eta^2)^2}}\,d\eta. \]
In case $t:-\pi/6\longmapsto 0$, by the similar arguments, we obtain 
\[\dfrac{z^5+z}{w^2}\,dz=i\,
\dfrac{-1}{\sqrt[3]{4\eta^2(3-4\eta^2)^2}}\,d\eta. \]
Therefore, we obtain 
\begin{align*}
\int_{A_1} \dfrac{z^5+z}{w^2}\,dz&=
\int^{\frac{\sqrt{3}}{2}}_0 
i\,\dfrac{e^{\frac{\pi}{3}\,i}}{\sqrt[3]{4\eta^2(3-4\eta^2)^2}}\,d\eta
+\int_{\frac{\sqrt{3}}{2}}^0 
i\,\dfrac{-1}
{\sqrt[3]{4\eta^2(3-4\eta^2)^2}}\,d\eta\\
&= i\,(1+e^{\frac{\pi}{3}\,i})\,\int^{\frac{\sqrt{3}}{2}}_0\dfrac{d\eta}
{\sqrt[3]{4\eta^2(3-4\eta^2)^2}}\\
&=\dfrac{i\,(1+e^{\frac{\pi}{3}\,i})}{2\sqrt{3}}\,
\int_0^1\dfrac{dt}{\sqrt[3]{t^2(1-t^2)^2}}\quad 
\left(\eta=\dfrac{\sqrt{3}}{2}\,t\right)\\
&=\dfrac{i\,(1+e^{\frac{\pi}{3}\,i})}{4\sqrt{3}}\, 
\int^1_0 s^{-5/6}(1-s)^{-2/3}\,ds\quad (s=t^2)\\
&=\dfrac{i\,(1+e^{\frac{\pi}{3}\,i})}{4\sqrt{3}}\, B(1/3,1/6). 
\end{align*}
Similarly, $A_2$-period is calculated as follows: 
\begin{align*}
\int_{A_2} \dfrac{z^5+z}{w^2}\,dz&=
\int^{\frac{\sqrt{3}}{2}}_0 
i\,\dfrac{e^{\frac{\pi}{3}\,i}}{\sqrt[3]{4\eta^2(3-4\eta^2)^2}}\,d\eta
+\int_{\frac{\sqrt{3}}{2}}^0 
i\,\dfrac{-e^{\frac{2\pi}{3}\,i}}{\sqrt[3]{4\eta^2(3-4\eta^2)^2}}\,d\eta\\
&= i\,(e^{\frac{\pi}{3}\,i}+e^{\frac{2\pi}{3}\,i})\,\int^{\frac{\sqrt{3}}{2}}_0\dfrac{d\eta}
{\sqrt[3]{4\eta^2(3-4\eta^2)^2}}\\
&=\dfrac{i\,(e^{\frac{\pi}{3}\,i}+e^{\frac{2\pi}{3}\,i})}{4\sqrt{3}}\,B(1/3,1/6). 
\end{align*}

\begin{firth_{10}}
\end{firth_{10}}

First, we consider 
$A_1$-period. In case $t: 0\longmapsto\pi/6$, we set 
\[\eta=\dfrac{1}{2}\,\left(z^2+\dfrac{1}{z^2}\right)=
\cos (2t).\] 
Then, $d\eta=-\dfrac{1-z^4}{z^3}\,dz$, and thereby 
\[\dfrac{i\,(z^5-z)}{w^2}\,dz=\dfrac{-i\,z(1-z^4)}{w^2}\,\dfrac{-z^3}{1-z^4}\,d\eta
=i\,\dfrac{z^4}{w^2}\,d\eta.\] 
To calculate $z^4/w^2$, we consider 
\[\left(\dfrac{z^4}{w^2}\right)^3=\dfrac{z^{12}}{(z^{12}-1)^2}
=\dfrac{1}{z^{12}+\dfrac{1}{z^{12}}-2}.\]
Note that 
$z^2+\dfrac{1}{z^2}=2\eta$, 
$z^4+\dfrac{1}{z^4}=4\eta^2-2$, and thus 
\[z^6+\dfrac{1}{z^6}=(z^2+\dfrac{1}{z^2})(z^4-1+\dfrac{1}{z^4})
=2\eta(4\eta^2-3).\] 
It follows that  
\[z^{12}+\dfrac{1}{z^{12}}-2=(2\eta(4\eta^2-3))^2-4=
-4\,(1-\eta^2)(4\eta^2-1)^2.\]
So, we obtain 
\[\left(\dfrac{z^4}{w^2}\right)^3=-\dfrac{1}{4\,(1-\eta^2)(4\eta^2-1)^2}<0. \]
Here, we take a suitable branch for $z^4/w^2$. 
\[\dfrac{z^4}{w^2}\left(\dfrac{\pi}{12}\right)=\dfrac{e^{\frac{\pi}{3}\,i}}
{\sqrt[3]{4}}. \]
Choosing a branch $\sqrt[3]{ 4\,(1-\eta^2)(4\eta^2-1)^2}>0$, we get 
\[\dfrac{z^4}{w^2}=\dfrac{e^{\frac{\pi}{3}\,i}}{\sqrt[3]{ 4\,(1-\eta^2)(4\eta^2-1)^2}}.\] 
It follows that 
\[\dfrac{i\,(z^5-z)}{w^2}\,dz=i\,
\dfrac{e^{\frac{\pi}{3}\,i}}{\sqrt[3]{4\,(1-\eta^2)(4\eta^2-1)^2}}\,d\eta. \]
In case $t:-\pi/6\longmapsto 0$, by the similar arguments, we obtain 
\[\dfrac{i\,(z^5-z)}{w^2}\,dz=i\,
\dfrac{-1}{\sqrt[3]{4\,(1-\eta^2)(4\eta^2-1)^2}}\,d\eta. \]
Therefore, we obtain
\begin{align*}
\int_{A_1} \dfrac{i\,(z^5-z)}{w^2}\,dz&=
\footnotesize{\int^{\frac{1}{2}}_1 
i\,\dfrac{e^{\frac{\pi}{3}\,i}}{\sqrt[3]{4\,(1-\eta^2)(4\eta^2-1)^2}}\,d\eta
+\int_{\frac{1}{2}}^1 
i\,\dfrac{-1}{\sqrt[3]{4\,(1-\eta^2)(4\eta^2-1)^2}}\,d\eta}\\
&= -i\,(1+e^{\frac{\pi}{3}\,i})\,\int^1_{\frac{1}{2}}\dfrac{dt}
{\sqrt[3]{4\,(1-t^2)(4t^2-1)^2}}.
\end{align*}

Similarly, $A_2$-period is calculated as follows: 
\begin{align*}
\int_{A_2} \dfrac{i\,(z^5-z)}{w^2}\,dz&=
\footnotesize{\int^{\frac{1}{2}}_1 
i\,\dfrac{e^{\frac{\pi}{3}\,i}}{\sqrt[3]{4\,(1-\eta^2)(4\eta^2-1)^2}}\,d\eta
+\int_{\frac{1}{2}}^1 
i\,\dfrac{-e^{\frac{2\pi}{3}\,i}}{\sqrt[3]{4\,(1-\eta^2)(4\eta^2-1)^2}}\,d\eta}\\
&= -i\,(e^{\frac{\pi}{3}\,i}+e^{\frac{2\pi}{3}\,i})\,\int^1_{\frac{1}{2}}\dfrac{dt}
{\sqrt[3]{4\,(1-t^2)(4t^2-1)^2}}.
\end{align*}

\end{proof}

Here, we set 
\[\alpha=\dfrac{1}{6\sqrt[3]{2}}\,B(2/3,1/6),\>\>
\beta=\dfrac{1}{4\sqrt{3}}\,B(1/3,1/6),\>\>
\gamma=\int^{1}_{\frac{1}{2}}\dfrac{dt}{\sqrt[3]{4(1-t^2)(4t^2-1)^2}},\]
and determine a relation between $\beta$ and $\gamma$: 

\begin{lem}

$\beta=\sqrt{3}\,\gamma$.
\end{lem}

\begin{proof}

We consider the period of $\dfrac{z^5+z}{w^2}\,dz$ along $B_5=\varphi^5(A_1)$. 
The followings are the similar arguments as in Lemma 3.1. 
In case $t:5\pi/6\longmapsto \pi$, we set $\eta=\sin(2t)$. Note that 
\[\dfrac{z^4}{w^2}\left(\dfrac{11}{12}\,\pi\right)=\dfrac{e^{\frac{11\pi}{3}\,i}}
{\sqrt[3]{4}}=-\dfrac{e^{\frac{2\pi}{3}\,i}}{\sqrt[3]{4}}.\]
Choosing a branch $\sqrt[3]{4\eta^2(3-4\eta^2)^2}>0$, we get 
\[\dfrac{z^5+z}{w^2}\,dz=i\,\dfrac{-e^{\frac{2\pi}{3}\,i}}
{\sqrt[3]{4\eta^2(3-4\eta^2)^2}}\,d\eta. \]
Similarly, in case $t:-\pi\longmapsto -5\pi/6$, we obtain
\[\dfrac{z^5+z}{w^2}\,dz=i\,\dfrac{e^{\frac{\pi}{3}\,i}}
{\sqrt[3]{4\eta^2(3-4\eta^2)^2}}\,d\eta. \]
It follows that 
\begin{align}
\int_{B_5} \dfrac{z^5+z}{w^2}\,dz&=\int^{0}_{-\frac{\sqrt{3}}{2}} 
i\,\dfrac{-e^{\frac{2\pi}{3}\,i}}{\sqrt[3]{4\eta^2(3-4\eta^2)^2}}
\,d\eta+\int_{0}^{-\frac{\sqrt{3}}{2}} 
i\,\dfrac{e^{\frac{\pi}{3}\,i}}
{\sqrt[3]{4\eta^2(3-4\eta^2)^2}}\,d\eta\\
\nonumber &=-i\,(e^{\frac{2\pi}{3}\,i}+e^{\frac{\pi}{3}\,i}) \,\beta.
\end{align}
On the other hand, using the action of $\varphi$, we obtain 
\[(\varphi^*)^5=e^{\frac{10\pi}{3}\,i}\begin{pmatrix}
R\left(\dfrac{\pi}{2}\right)&0\\
0&R\left(-\dfrac{\pi}{3}\right)
\end{pmatrix}^5
=-e^{\frac{\pi}{3}\,i}\begin{pmatrix}
0&-1&0&0\\
1&0&0&0\\
0&0&\dfrac{1}{2}&-\dfrac{\sqrt{3}}{2}\\
0&0&\dfrac{\sqrt{3}}{2}&\dfrac{1}{2}
\end{pmatrix}.
\]
Hence, 
\begin{align}
\nonumber \int_{B_5} \dfrac{z^5+z}{w^2}\,dz&=\int_{\varphi^5(A_1)} \dfrac{z^5+z}{w^2}\,dz
=\int_{A_1}(\varphi^*)^5 \left(\dfrac{z^5+z}{w^2}\,dz\right)\\
&\nonumber =-e^{\frac{\pi}{3}\,i}
\begin{pmatrix}
\dfrac{1}{2}&-\dfrac{\sqrt{3}}{2}
\end{pmatrix}
\begin{pmatrix}
\int_{A_1} \dfrac{z^5+z}{w^2}\,dz\\
\int_{A_1} \dfrac{i\,(z^5-z)}{w^2}\,dz
\end{pmatrix}\\
\nonumber &=-e^{\frac{\pi}{3}\,i}
\begin{pmatrix}
\dfrac{1}{2}&-\dfrac{\sqrt{3}}{2}
\end{pmatrix}
\begin{pmatrix}
i\,(1+e^{\frac{\pi}{3}\,i})\,\beta\\
-i\,(1+e^{\frac{\pi}{3}\,i})\,\gamma\\
\end{pmatrix}\\
\nonumber &=-e^{\frac{\pi}{3}\,i}\left(\dfrac{i}{2}(1+e^{\frac{\pi}{3}\,i})\,\beta
+\dfrac{\sqrt{3}}{2}\,i\,(1+e^{\frac{\pi}{3}\,i})\,\gamma\right)\\
&=-i\,(e^{\frac{\pi}{3}\,i}+e^{\frac{2\pi}{3}\,i})\left(\dfrac{\beta}{2}+
\dfrac{\sqrt{3}}{2}\,\gamma\right).
\end{align}
$(5)=(6)$ implies $\beta=\sqrt{3}\,\gamma$ and this completes the proof. 
\end{proof}
By Lemma 3.1, 3.2, and the action of $\varphi$, the period matrix $\Re\Omega$ of $f$ is given 
as follows: 
\[\Re\Omega=\Re\begin{pmatrix}
\Omega_1,\>\Omega_2,\>\Omega_3,\>\Omega_4,\>\Omega_5,\>\Omega_6,\>\Omega_7
\end{pmatrix},\]
where
\[\Omega_1=
\begin{pmatrix}
(1+e^{\frac{\pi}{3}\,i})\,\begin{pmatrix}
\alpha\\           
-\alpha\\
\sqrt{3}\,i\,\gamma\\
-i\,\gamma
\end{pmatrix},
(e^{\frac{2\pi}{3}\,i}+e^{\frac{\pi}{3}\,i})\,\begin{pmatrix}
\alpha\\           
-\alpha\\
\sqrt{3}\,i\,\gamma\\
-i\,\gamma
\end{pmatrix},
(-1+e^{\frac{2\pi}{3}\,i})\,\begin{pmatrix}
\alpha\\           
\alpha\\
0\\
-2\,i\,\gamma
\end{pmatrix}
\end{pmatrix},\]
\[\Omega_2=
\begin{pmatrix}
-(1+e^{\frac{\pi}{3}\,i})\,\begin{pmatrix}
\alpha\\           
\alpha\\
0\\
-2\,i\,\gamma
\end{pmatrix},
(e^{\frac{2\pi}{3}\,i}+e^{\frac{\pi}{3}\,i})\,\begin{pmatrix}
\alpha\\           
-\alpha\\
\sqrt{3}\,i\,\gamma\\
i\,\gamma
\end{pmatrix},
(-1+e^{\frac{2\pi}{3}\,i})\,\begin{pmatrix}
\alpha\\           
-\alpha\\
\sqrt{3}\,i\,\gamma\\
i\,\gamma
\end{pmatrix}
\end{pmatrix},\]
\[\Omega_3=
\begin{pmatrix}
-(1+e^{\frac{\pi}{3}\,i})\,\begin{pmatrix}
\alpha\\           
\alpha\\
\sqrt{3}\,i\,\gamma\\
-i\,\gamma
\end{pmatrix},
-(e^{\frac{2\pi}{3}\,i}+e^{\frac{\pi}{3}\,i})\,\begin{pmatrix}
\alpha\\           
\alpha\\
\sqrt{3}\,i\,\gamma\\
-i\,\gamma
\end{pmatrix},
(-1+e^{\frac{2\pi}{3}\,i})\,\begin{pmatrix}
\alpha\\           
-\alpha\\
0\\
2\,i\,\gamma
\end{pmatrix}
\end{pmatrix},\]
\[\Omega_4=
\begin{pmatrix}
-(1+e^{\frac{\pi}{3}\,i})\,\begin{pmatrix}
\alpha\\           
-\alpha\\
0\\
2\,i\,\gamma
\end{pmatrix},
-(e^{\frac{2\pi}{3}\,i}+e^{\frac{\pi}{3}\,i})\,\begin{pmatrix}
\alpha\\           
\alpha\\
\sqrt{3}\,i\,\gamma\\
i\,\gamma
\end{pmatrix},
(1-e^{\frac{2\pi}{3}\,i})\,\begin{pmatrix}
\alpha\\           
\alpha\\
\sqrt{3}\,i\,\gamma\\
i\,\gamma
\end{pmatrix}
\end{pmatrix},\]
\[\Omega_5=
\begin{pmatrix}
-(1+e^{\frac{\pi}{3}\,i})\,\begin{pmatrix}
\alpha\\           
-\alpha\\
-\sqrt{3}\,i\,\gamma\\
i\,\gamma
\end{pmatrix},
-(e^{\frac{2\pi}{3}\,i}+e^{\frac{\pi}{3}\,i})\,\begin{pmatrix}
\alpha\\           
-\alpha\\
-\sqrt{3}\,i\,\gamma\\
i\,\gamma
\end{pmatrix},
(1-e^{\frac{2\pi}{3}\,i})\,\begin{pmatrix}
\alpha\\           
\alpha\\
0\\
2\,i\,\gamma
\end{pmatrix}
\end{pmatrix},\]
\[\Omega_6=
\begin{pmatrix}
(1+e^{\frac{\pi}{3}\,i})\,\begin{pmatrix}
\alpha\\           
\alpha\\
0\\
2\,i\,\gamma
\end{pmatrix},
-(e^{\frac{2\pi}{3}\,i}+e^{\frac{\pi}{3}\,i})\,\begin{pmatrix}
\alpha\\           
-\alpha\\
-\sqrt{3}\,i\,\gamma\\
-i\,\gamma
\end{pmatrix},
(1-e^{\frac{2\pi}{3}\,i})\,\begin{pmatrix}
\alpha\\           
-\alpha\\
-\sqrt{3}\,i\,\gamma\\
-i\,\gamma
\end{pmatrix}
\end{pmatrix},\]
\[\Omega_7=\begin{pmatrix}
(1+e^{\frac{\pi}{3}\,i})\,\begin{pmatrix}
\alpha\\           
\alpha\\
-\sqrt{3}\,i\,\gamma\\
i\,\gamma
\end{pmatrix},
(e^{\frac{2\pi}{3}\,i}+e^{\frac{\pi}{3}\,i})\,\begin{pmatrix}
\alpha\\           
\alpha\\
-\sqrt{3}\,i\,\gamma\\
i\,\gamma
\end{pmatrix}
\end{pmatrix}.\]
We now review the lattice transformation. 
A lattice $\Lambda$ in a real vector space $\textbf{R}^n$ is a 
discrete subgroup of maximal rank in $\textbf{R}^n$. Let 
$\{u_1,u_2,\cdots,u_m\}$ $(m\geq n)$ be a sequence of vectors 
which span $\textbf{R}^n$. In general, $\{u_1,u_2,\cdots,u_m\}$ 
are not lattice vectors. 
\begin{prop}[section 6 in \cite{E}]$\,$

$\{u_1,u_2,\cdots,u_m\}$ are lattice vectors if and only if there 
exist lattice vectors $\{v_1,v_2,\cdots,v_n\}$ such that 
\begin{align*}
&\{v_1,v_2,\cdots,v_n\}=\{u_1,u_2,\cdots,u_m\}\,G_1,\\
&\{u_1,u_2,\cdots,u_m\}=\{v_1,v_2,\cdots,v_n\}\,G_2,
\end{align*}
where $G_1$ is an $(m,n)$-matrix and $G_2$ is an $(n,m)$-matrix 
whose components are integers. 
\end{prop}
First, we set the following four matrices: 
\[\Omega_8=\begin{pmatrix}
-2\,(e^{\frac{2\,\pi}{3}\,i}-1)\,\alpha &2\,(1+e^{\frac{\pi}{3}\,i})\,\alpha &
(e^{\frac{\,\pi}{3}\,i}+e^{\frac{2\,\pi}{3}\,i})\,\alpha&
-(1+e^{\frac{\,\pi}{3}\,i})\,\alpha\\
0&0&-(e^{\frac{\,\pi}{3}\,i}+e^{\frac{2\,\pi}{3}\,i})\,\alpha&
(1+e^{\frac{\,\pi}{3}\,i})\,\alpha\\
0&0&0&0\\
0&0&0&0
\end{pmatrix},\]
\[\Omega_9=\footnotesize{\begin{pmatrix}
0&0&0&0\\
0&0&0&0\\
-2\sqrt{3}\,i\,(e^{\frac{\,\pi}{3}\,i}+e^{\frac{2\,\pi}{3}\,i})\,\gamma &0&
\sqrt{3}\,i\,(1+e^{\frac{\,\pi}{3}\,i})\gamma&
\sqrt{3}\,i\,(e^{\frac{\,\pi}{3}\,i}+e^{\frac{2\,\pi}{3}\,i})\,\gamma\\
0&2\,i\,(1+e^{\frac{\,\pi}{3}\,i})\,\gamma&
i\,(1+e^{\frac{\pi}{3}\,i})\gamma
&-i\,(e^{\frac{\,\pi}{3}\,i}+e^{\frac{2\,\pi}{3}\,i})\,\gamma
\end{pmatrix}},\]
\[G^{\Omega}_{1}=\tiny{\begin{pmatrix}
0& 0& 0& 0& 0& 0& 0& 0\\
0& 0& 1& -1& -1& -1& 0& 1\\
-1& 0& 0& 1& 0& 1& 0& -1\\
0& -1& 0& -1& 0& -1& 0& 1\\
0& 0& -1& 1& 0& 1& 0& -1\\
0& 0& 0& 1& 0& 1& 0& 0\\
0& 0& 0& 0& 0& 0& 0& 0\\
0& 0& 0& 1& 0& 1& 0& -1\\
-1& 0& 1& -1& 0& -1& 0& 0\\
0& -1& -1& 0& 0& -1& -1& 0\\
0& 0& 0& 0& 0& 0& 0& 0\\
0& 0& 0& 0& 0& 0& 0& 0\\
0& 0& 0& 1& 0& 1& 1& 0\\
0& 0& 0& 0& 0& 0& 0& 0\\
0& 0& 0& 0& 0& 0& 0& 0\\
0& 0& 0& 0& 0& 0& 0& 0\\
0& 0& 0& 0& -1& 0& 0& 0\\
0& 0& 0& 0& 0& 0& 0& 0\\
0& 0& 0& 0& 0& 0& 0& 0\\
0& 0& 0& 0& 0& 0& 0& 0
\end{pmatrix}},
(G^{\Omega}_{2})^T=\tiny{\begin{pmatrix}
0& 0& 0& -1& 0& -1& 1& 0\\
0& 0& 1& 0& 0& 0& 0& 1\\
-1& 0& -1& -1& 1& 1& 0& 2\\
0& -1& 0& -1& 0& 1& 0& 0\\
0& 0& 1& 0& -1& 0& 0& -1\\
0& 0& 1& 1& -1& 0& 1& -1\\
0& -1& 0& -1& 0& 1& -1& 0\\
1& -1& 1& 0& 0& 0& 0& -1\\
0& 0& 1& 1& -1& -1& 0& -2\\
0& 0& 0& 1& 0& -1& 0& 0\\
1& -1& 1& 0& 1& 0& 0& 1\\
1& 0& 1& 1& 1& 0& 1& 1\\
0& 0& 0& 1& 0& -1& 1& 0\\
0& 0& -1& 0& 0& 0& 0& 1\\
1& 0& 1& 1& 1& 1& 0& 2\\
0& 1& 0& 1& 0& 1& 0& 0\\
0& 0& -1& 0& -1& 0& 0& -1\\
0& 0& -1& -1& -1& 0& -1& -1\\
0& 1& 0& 1& 0& 1& -1& 0\\
-1& 1& -1& 0& 0& 0& 0& -1
\end{pmatrix}}.\]
Then, 
\[\Omega \,G^{\Omega}_1 = \begin{pmatrix}
\Omega_8,\Omega_9
\end{pmatrix},\quad \begin{pmatrix}
\Omega_8,\Omega_9
\end{pmatrix}G^{\Omega}_2=\Omega\]
hold. 
 
Taking the real and imaginary parts respectively, we get 
\begin{align*}
(\Omega_8,\Omega_9)&=\underbrace{\begin{pmatrix}
3\,\alpha& 3\,\alpha & 0 & -\dfrac{3}{2}\,\alpha & 0 & 0& 0 & 0\\
0& 0 & 0 & \dfrac{3}{2}\,\alpha & 0 & 0& 0 &0\\
0 & 0 & 0 & 0 & 6\,\gamma& 0 & -\dfrac{3}{2}\,\gamma & -3\,\gamma\\
0 & 0 & 0 & 0 & 0 & -\sqrt{3}\,\gamma & -\dfrac{\sqrt{3}}{2}\,\gamma & \sqrt{3}\,\gamma
\end{pmatrix}}_{\Omega_R}\\
&\quad+\,i
\underbrace{\begin{pmatrix}
-\sqrt{3}\,\alpha & \sqrt{3}\,\alpha & \sqrt{3}\,\alpha & 
-\dfrac{\sqrt{3}}{2}\,\alpha& 0 & 0 &0 & 0 \\
0 & 0 & -\sqrt{3}\,\alpha & \dfrac{\sqrt{3}}{2}\,\alpha& 0 & 0 & 0 & 0\\
0 & 0 & 0 & 0 & 0 & 0 & \dfrac{3\,\sqrt{3}}{2}\,\gamma& 0 \\
0 & 0 & 0 & 0 & 0 & 3\,\gamma & \dfrac{3}{2}\,\gamma& 0 
\end{pmatrix}}_{\Omega_I}.
\end{align*}

Next, we consider the followings: 
\[G^R_1=\begin{pmatrix}
1&1&0&0\\
0&0&0&0\\
0&0&0&0\\
0&1&0&0\\
0&0&0&0\\
0&0&1&0\\
0&0&-2&-1\\
0&0&0&0\\
\end{pmatrix},\quad
G^R_2 =\begin{pmatrix}
1&1&0&-1&0&0&0&0\\
0&0&0&1&0&0&0&0\\
0&0&0&0&2&1&0&-2\\
0&0&0&0&0&-2&-1&2
\end{pmatrix},\]
\[G^I_1=\begin{pmatrix}
-1&-1&0&0\\
0&0&0&0\\
0&0&0&0\\
0&1&0&0\\
0&0&0&0\\
0&0&-1&0\\
0&0&2&1\\
0&0&0&0
\end{pmatrix},\quad
G^I_2=\begin{pmatrix}
-1&1&2&-1&0&0&0&0\\
0&0&-2&1&0&0&0&0\\
0&0&0&0&0&-1&0&0\\
0&0&0&0&0&2&1&0
\end{pmatrix},\]
\[\Lambda=\begin{pmatrix}
3\,\alpha & \dfrac{3}{2}\,\alpha & 0 & 0\\
0 & \dfrac{3}{2}\,\alpha & 0 & 0\\
0 & 0 & 3\,\gamma & \dfrac{3}{2}\,\gamma\\
0 & 0 & 0 & \dfrac{\sqrt{3}}{2}\,\gamma
\end{pmatrix},
\quad
\Lambda_{\pi/2}=\begin{pmatrix}
\sqrt{3}\,\alpha & \dfrac{\sqrt{3}}{2}\,\alpha & 0 & 0\\
0 & \dfrac{\sqrt{3}}{2}\,\alpha & 0 & 0\\
0 & 0 & 3\sqrt{3}\,\gamma & \dfrac{3\sqrt{3}}{2}\,\gamma\\
0 & 0 & 0 & \dfrac{3}{2}\,\gamma
\end{pmatrix}.
\]
Then, we obtain 
\[\Omega_R\,G^R_1=\Lambda,\quad \Lambda\,G^R_2=\Omega_R,\quad 
\Omega_I\,G^I_1=\Lambda_{\pi/2},\quad \Lambda_{\pi/2}\,G^I_2=\Omega_I.\]
So, by Proposition 3.1, $\Omega_R$ (resp. $\Omega_I$) determines a lattice given by 
$\Lambda$ (resp. $\Lambda_{\pi/2}$). Therefore, we can define a minimal surface 
$f:M_{10}\longrightarrow \textbf{R}^4/\Lambda$ and the conjugate surface 
$f_{\pi/2}:M_{10}\longrightarrow \textbf{R}^4/\Lambda_{\pi/2}$. 

Finally, we consider the associate surfaces $f_{\theta}$ of $f$. 
The period matrix of $f_{\theta}$ is given by 
\begin{align*}
\Re\{e^{i\,\theta}\,\Omega\}&=
\cos\theta\begin{pmatrix}
\Omega_{10} & \Omega_{11} 
\end{pmatrix},
\end{align*}
where
\begin{align*}
\Omega_{10}&=\footnotesize{\begin{pmatrix}
(3+\sqrt{3}\,\tan\theta)\,\alpha & (3-\sqrt{3}\,\tan\theta)\,\alpha & 
-\sqrt{3}\,\alpha \,\tan\theta& -\dfrac{\alpha}{2}\,(3-\sqrt{3}\,\tan\theta) \\ 
0 & 0 & \sqrt{3}\,\alpha \,\tan\theta& \dfrac{\alpha}{2}\,(3-\sqrt{3}\,\tan\theta) \\ 
0 & 0 & 0 & 0 \\
0 & 0 & 0 & 0  
\end{pmatrix}}, \\
\Omega_{11}&=\begin{pmatrix}
0 & 0 & 0 & 0\\
0 & 0 & 0 & 0\\
6\,\gamma & 0 & -\dfrac{3\,\gamma}{2}\,(1+\sqrt{3}\,\tan\theta) & -3\,\gamma\\
0&-\sqrt{3}\,\gamma\,(1+\sqrt{3}\,\tan\theta) &
-\dfrac{\sqrt{3}\,\gamma}{2}\,(1+\sqrt{3}\,\tan\theta) & \sqrt{3}\,\gamma
\end{pmatrix}.
\end{align*}
Setting $\sqrt{3}\,\tan\theta=\dfrac{m}{n}\in \textbf{Q}$, we can show that 
${\rm rank}_{\textbf{Q}}\Re\{e^{i\,\theta}\,\Omega\}=4$. 
This implies (3) in Theorem 1.1, and thereby each $f_{\theta}$ is well-defined. 
Moreover, the associate surfaces are well-defined for a countable dense set of angles 
$e^{i\,\theta}\in S^1$ because $\theta$ is parametrized by the rational number. 

\subsection{Homological triviality}

\begin{lem}
$f(M_{10})$ is homologous to $0$ in $\textbf{R}^4/\Lambda$. 
\end{lem}
\begin{proof}
Setting 
\begin{align*}
(x^1,\,x^2,\,x^3,\,x^4)=\Re \int_{p_0}^{p} \left(
\dfrac{1-z^6}{w^2},\,
\dfrac{i\,(1+z^6)}{w^2},\,
\dfrac{z^5+z}{w^2},\,
\dfrac{i\,(z^5-z)}{w^2}
\right)^T\,dz,
\end{align*} 
we obtain 
\begin{align*}
\begin{pmatrix}
dx^1\\
dx^2\\
dx^3\\
dx^4\end{pmatrix}
=\dfrac{1}{2}\,\begin{pmatrix}
\dfrac{1-z^6}{w^2}\,dz+\dfrac{1-\bar{z}^6}{\bar{w}^2}\,d\bar{z}\\
\dfrac{i\,(1+z^6)}{w^2}\,dz-\dfrac{i\,(1+\bar{z}^6)}{\bar{w}^2}\,d\bar{z}\\
\dfrac{z^5+z}{w^2}\,dz+\dfrac{\bar{z}^5+\bar{z}}{\bar{w}^2}\,d\bar{z}\\
\dfrac{i\,(z^5-z)}{w^2}\,dz-\dfrac{i\,(\bar{z}^5-\bar{z})}{\bar{w}^2}\,d\bar{z}\end{pmatrix}.
\end{align*}

It follows that 
\begin{align*}
dx^1\wedge dx^2&=
\dfrac{1}{4}\,\left(\dfrac{1-z^6}{w^2}\,dz+\dfrac{1-\bar{z}^6}{\bar{w}^2}\,d\bar{z}\right)
\wedge\left(\dfrac{i\,(1+z^6)}{w^2}\,dz-\dfrac{i\,(1+\bar{z}^6)}{\bar{w}^2}\,d\bar{z}\right)\\
&=-\dfrac{i}{2}\,\dfrac{1-|z|^{12}}{|w|^4}\,dz\wedge d\bar{z},\\
dx^1\wedge dx^3&=
\dfrac{1}{4}\,\left(\dfrac{1-z^6}{w^2}\,dz+\dfrac{1-\bar{z}^6}{\bar{w}^2}\,d\bar{z}\right)
\wedge\left(\dfrac{z^5+z}{w^2}\,dz+\dfrac{\bar{z}^5+\bar{z}}{\bar{w}^2}\,d\bar{z}\right)\\
&=\dfrac{1}{4}\,\dfrac{-(z-\bar{z})(1+|z|^{10})-(z^5-\bar{z}^5)(1+|z|^2)}{|w|^4}\,
dz\wedge d\bar{z},\\
dx^1\wedge dx^4&=
\dfrac{1}{4}\,\left(\dfrac{1-z^6}{w^2}\,dz+\dfrac{1-\bar{z}^6}{\bar{w}^2}\,d\bar{z}\right)
\wedge\left(\dfrac{i\,(z^5-z)}{w^2}\,dz-\dfrac{i\,(\bar{z}^5-\bar{z})}{\bar{w}^2}\,d\bar{z}
\right)\\
&=-\dfrac{i}{4}\,\dfrac{-(z+\bar{z})(1+|z|^{10})+(z^5+\bar{z}^5)(1+|z|^2)}{|w|^4}\,
dz\wedge d\bar{z},\\
dx^2\wedge dx^3&=
\dfrac{1}{4}\,\left(\dfrac{i\,(1+z^6)}{w^2}\,dz-\dfrac{i\,(1+\bar{z}^6)}{\bar{w}^2}\,
d\bar{z}\right)\wedge
\left(\dfrac{z^5+z}{w^2}\,dz+\dfrac{\bar{z}^5+\bar{z}}{\bar{w}^2}\,d\bar{z}\right)\\
&=\dfrac{i}{4}\,\dfrac{(z+\bar{z})(1+|z|^{10})+(z^5+\bar{z}^5)(1+|z|^2)}{|w|^4}\,
dz\wedge d\bar{z},\\
dx^2\wedge dx^4&=
\dfrac{1}{4}\,\left(\dfrac{i\,(1+z^6)}{w^2}\,dz-\dfrac{i\,(1+\bar{z}^6)}{\bar{w}^2}\,
d\bar{z}\right)\wedge
\left(\dfrac{i\,(z^5-z)}{w^2}\,dz-\dfrac{i\,(\bar{z}^5-\bar{z})}{\bar{w}^2}\,d\bar{z}\right)\\
&=\dfrac{1}{4}\,\dfrac{(z-\bar{z})(1+|z|^{10})-(z^5-\bar{z}^5)(1+|z|^2)}{|w|^4}\,
dz\wedge d\bar{z},\\
dx^3\wedge dx^4&=
\dfrac{1}{4}\,\left(\dfrac{z^5+z}{w^2}\,dz+\dfrac{\bar{z}^5+\bar{z}}{\bar{w}^2}\,d\bar{z}
\right)
\wedge\left(\dfrac{i\,(z^5-z)}{w^2}\,dz-\dfrac{i\,(\bar{z}^5-\bar{z})}{\bar{w}^2}\,d\bar{z}
\right)\\
&=-\dfrac{i}{2}\,\dfrac{-|z|^2+|z|^{10}}{|w|^4}\,dz\wedge d\bar{z}.
\end{align*}
Setting $z=re^{i\,\theta}\,\,(0\leq r \leq \infty,\,0\leq\theta\leq2\pi)$, we get 
$dz\wedge d\bar{z}=-2\,i\,rdr\wedge d\theta$. 
First, 
\begin{align*}
\int_{M_{10}}dx^1\wedge dx^2&=3\int_{\theta=0}^{2\pi}\int_{r=0}^{\infty} 
\dfrac{i}{2}\,\dfrac{1-r^{12}}
{\sqrt[3]{(r^{24}-2r^{12}\,\cos(12\,\theta)+1)^2}}\,2ir\,dr\,d\theta\\
&=-3\int_{\theta=0}^{2\pi}\left(\int_{r=0}^{1} +\int_{r=1}^{\infty}\right)\\
&=-3\int_{\theta=0}^{2\pi}\left(\int_{r=0}^{1} -\int_{r'=0}^{1}\right) \qquad
(r'=1/r)\\
&=0.
\end{align*}
Similarly, we can see that $\int_{M_{10}}dx^3\wedge dx^4=0$. 
Next, we obtain 
\begin{align*}
\int_{M_{10}}dx^1\wedge dx^3
&=-3\int_{\theta=0}^{2\pi}\int^{\infty}_{r=0}
\dfrac{r\sin\theta(1+r^{10})+r^5\,\sin(5\,\theta)(1+r^2)}
{\sqrt[3]{(r^{24}-2r^{12}\,\cos(12\,\theta)+1)^2}}\,rdr\,d\theta
\end{align*}
\begin{align*}
&=-3\left(\int_{\theta=0}^{\pi}\int^{\infty}_{r=0}
+\int_{\theta=\pi}^{2\pi}\int^{\infty}_{r=0}\right) \\
&=-3\left(\int_{\theta=0}^{\pi}\int^{\infty}_{r=0}-\int_{\theta'=0}^{\pi}\int^{\infty}_{r=0} \right)
\qquad(\theta'=\theta-\pi)\\
&=0.
\end{align*}
Similarly, 
$ \int_{M_{10}}dx^1\wedge dx^4 
=\int_{M_{10}}dx^2\wedge dx^3 =\int_{M_{10}}dx^2 \wedge dx^4 =0$ follows. 
Therefore, $[f(M_{10})]=0$ holds. 
\end{proof}

\subsection{Symmetry}

\begin{lem}
$f$ has only reducible symmetry, and moreover, the 
maximal symmetry of $f$ is the dihedral group $D_{12}$. 
\end{lem}
\begin{proof}
Let $S_f(M_{10})$ be an arbitrary symmetry and, for every $\phi\in S_f(M_{10})$, 
the following diagram commutes: 

\begin{picture}(300,40)
 \put(60,30){\vector(1,0){190}}
 \put(30,27){$M_{10}$}
 \put(260,27){$\textbf{R}^4/\Lambda$}
 \put(60,-25){\vector(1,0){190}}
 \put(30,-27){$M_{10}$}
 \put(36,18){\vector(0,-1){30}}
 \put(274,18){\vector(0,-1){30}}
 \put(153,0){$\circlearrowright$}
 \put(260,-27){$\textbf{R}^4/\Lambda$}
 \put(285,10){$\exists$ affine transformation}
 \put(20,0){$\phi$} \put(285,-3){$A\in O(4,\textbf{R}),\> \textbf{t}\in\textbf{R}^4$}
 \put(153,35){$f$}\put(153,-35){$f$}
\end{picture}

$\,$\\
\\
\\
Here, we introduce the automorphism $j$ defined by $j(z,w)=(z,e^{\frac{2\pi}{3}\,i}w)$. 
Because the Gauss map $G$ is $j$-invariant, $\phi$ induces an automorphism $\phi'$ of $S^2$:

\begin{picture}(300,40)
 \put(60,30){\vector(1,0){190}}
 \put(30,27){$M_{10}$}
 \put(260,27){$Q_2\subset \textbf{C}P^3$}
 \put(60,-115){\vector(1,0){190}}
 \put(30,-117){$M_{10}$}
 \put(36,18){\vector(0,-1){120}}
 \put(274,18){\vector(0,-1){120}}
 \put(260,-117){$Q_2\subset \textbf{C}P^3$}
 \put(20,-40){$\phi$}\put(285,-50){$A$}
 \put(148,35){$G$}\put(148,-125){$G$}
 \put(130,0){$M_{10}/j\cong S^2$}
 \put(130,-90){$M_{10}/j\cong S^2$}
 \put(60,20){\vector(4,-1){60}}
 \put(190,5){\vector(4,1){60}}
 \put(65,5){$/j$}\put(200,0){embedding}
 \put(60,-105){\vector(4,1){60}}
 \put(190,-90){\vector(4,-1){60}}
 \put(80,-90){$/j$}\put(200,-90){embedding}
 \put(155,-10){\vector(0,-1){65}} 
 \put(90,-40){$\circlearrowright$}
 \put(210,-40){$\circlearrowright$}
 \put(130,-40){$\exists$ $\phi'$}
\end{picture}

$\,$\\
\\
\\
\\
\\
\\
\\
\\
\\
\\
Let $S'_f(M_{10})$ be a subgroup of automorphims of $M_{10}$ 
which induce automorphisms of $M_{10}/j\cong S^2 \cong \textbf{C} \cup \{ \infty \}$. 
Then, we obtain $S_f(M_{10}) \subset S'_f(M_{10})$. 
Now let $p_\alpha=(e^{\frac{\pi}{6}\,\alpha\,i},0)$ 
($1\leq \alpha \leq 12$) be branch points of $/j$. Every element of 
$S'_f(M_{10})$ induces an automorphism of $\textbf{C} \cup \{ \infty \}$ which 
preserves $\{e^{\frac{\pi}{6}\,\alpha\,i}\}_{\alpha=1}^{12}$. 
It is easy to verify that a subgroup of automorphisms of $\textbf{C} \cup \{ \infty \}$ which 
preserve $\{e^{\frac{\pi}{6}\,\alpha\,i}\}_{\alpha=1}^{12}$ is generated by 
$z\longmapsto e^{\frac{\pi}{6}\,i}\,z$ and $z\longmapsto 1/z$. 
To lift these automorphisms, we consider the automorphism $\varphi'$ of $M_{10}$ given by 
$\varphi'(z,w)=(1/z,\,e^{\frac{\pi}{3}\,i}\,w/z^4)$. 
Then, $S'_f(M_{10})$ is generated by $j$, $\varphi$, and $\varphi'$. 
Here, setting 
$\phi_1=\varphi \circ j^2$ and $\phi_2(z,w)=\varphi' \circ j$, we obtain 
\begin{align*}
\phi_1^*f=\begin{pmatrix} 
R\left(\dfrac{\pi}{2}\right)&0\\
0&R\left(-\dfrac{\pi}{3}\right) \end{pmatrix}f, \quad
\phi_2^*f=\begin{pmatrix} 1&0&0&0\\
0&-1&0&0\\ 
0& 0& -1&0\\ 
0& 0&0 &1 \end{pmatrix}f.
\end{align*}
Obviously, $\phi_1$ and $\phi_2$ are 
reducible actions respectively and generate the dihedral group $D_{12}$. 
Note that $j$ does not induce any affine transformation of the torus. Hence, 
$S_f(M_{10})$ is a subgroup of $D_{12}$ generated by $\phi_1$ and $\phi_2$. 
Therefore, $S_f(M_{10})$ has only reducible symmetry and the maximal symmetry is $D_{12}$ 
generated by $\phi_1$ and $\phi_2$. 
\end{proof}
By the results in 3.2, Lemma 3.3, and 3.4, the proof of Main Theorem 2 is complete.

\providecommand{\bysame}{\leavevmode\hbox to3em{\hrulefill}\thinspace}
\providecommand{\MR}{\relax\ifhmode\unskip\space\fi MR }
\providecommand{\MRhref}[2]{%
  \href{http://www.ams.org/mathscinet-getitem?mr=#1}{#2}
}
\providecommand{\href}[2]{#2}

\end{document}